\documentclass[10pt,reqno]{amsart}
\usepackage{amssymb,amscd}   

\newtheorem{theorem}{Theorem}[section]
\newtheorem{proposition}[theorem]{Proposition} \newtheorem{corollary}[theorem]{Corollary}

\theoremstyle{remark}

\newtheorem{remarks}[theorem]{Remarks}

\theoremstyle{definition}

\begin{document}
\title[Peak Phenomena for Non-commutative $H^\infty$]{On peak phenomena for non-commutative $H^\infty$}
\author[Y.~Ueda]{Yoshimichi UEDA$^*$}
\address{
Graduate School of Mathematics, 
Kyushu University, 
Fukuoka, 810-8560, Japan
}
\email{ueda@math.kyushu-u.ac.jp}
\thanks{$^*\,$Supported in part by Grant-in-Aid for Young Scientists (B)17740096.}
\dedicatory{Dedicated to Professor Fumio Hiai on the occasion of his 60th birthday} 
\maketitle

\begin{abstract} A non-commutative extension of Amar and Lederer's peak set result is given. As its simple applications it is shown that any non-commutative $H^\infty$-algebra $H^\infty(M,\tau)$ has unique predual, and moreover some restriction in some of the results of Blecher and Labuschagne are removed, making them hold in full generality. 
\end{abstract}

\allowdisplaybreaks{

\section{Introduction} 

Let $H^\infty(\mathbb{D})$ be the Banach algebra of all bounded analytic functions on the unit disk $\mathbb{D}$ equipped with the supremum norm $\Vert\cdot\Vert_\infty$. It is known (but non-trivial) that $H^\infty(\mathbb{D})$ can be regarded as a closed subalgebra of $L^\infty(\mathbb{T})$ by $f(e^{\sqrt{-1}\theta}) := \lim_{r\nearrow1}f(re^{\sqrt{-1}\theta})$ a.e.~$\theta$. Then, $L^\infty(\mathbb{T})$ is isometrically isomorphic to $C(X)$ with a certain compact Hausdorff space $X$ via the Gel'fand representation $f \mapsto \hat{f}$, and the linear functional $f \in H^\infty(\mathbb{D}) \mapsto \frac{1}{2\pi}\int_0^{2\pi} f(e^{\sqrt{-1}\theta})\,d\theta$ is known to admit a unique representing measure $m$ on $X$ so that $\frac{1}{2\pi}\int_0^{2\pi} f(e^{\sqrt{-1}\theta})\,d\theta = \int_X \hat{f}(x)\,m(dx)$ holds. In this setup, Amar and Lederer \cite{Amar and Lederer:CRParis71} proved that any closed subset $F \subset X$ with $m(F)=0$ admits $f \in H^{\infty}(\mathbb{D})$ with $\Vert f\Vert_\infty \le 1$ such that $P := \{x \in X : \hat{f}(x) = 1\} = \{x \in X : |\hat{f}(x)| = 1 \}$ contains $F$ and $m(P) = 0$ still holds. This is a key in any existing proof of the uniqueness of predual of $H^\infty(\mathbb{D})$. The reader can find some information on Amar and Lederer's result in \cite[\S6]{Pelczynski:CBMS76} and also see \cite{Barbey:ArchMath75}.  

The main purpose of these notes is to provide an analogous fact of the above-mentioned result of Amar and Lederer for non-commutative $H^\infty$-algebras introduced by Arveson \cite{Arveson:AJM67} in the 60's under the name of finite maximal subdiagonal algebras. Here a non-commutative $H^\infty$-algebra means a $\sigma$-weakly closed (possibly non-self-adjoint) unital subalgebra $A$ of a finite von Neuamnn algebra $M$ with a faithful normal tracial state $\tau$ satisfying the following conditions: 
\begin{itemize} 
\item the unique $\tau$-preserving (i.e., $\tau\circ E = \tau$) conditional expectation $E : M \rightarrow D := A\cap A^*$ is multiplicative on $A$; 
\item the $\sigma$-weak closure of $A+A^*$ is exactly $M$, 
\end{itemize}
where $A^* := \{a^* \in M : a \in A\}$. (Remark here that an important work due to Exel \cite{Exel:AJM88} plays an important r\^{o}le behind this simple definition.) In what follows we write $A = H^\infty(M,\tau)$ and call $D$ the diagonal subalgebra. Recently, in their series of papers Blecher and Labuschagne established many fundamental properties of these non-commutative $H^\infty$-algebras, analogous to classical theories modeled after $H^\infty(\mathbb{D})$, all of which are nicely summarized in \cite{BlecherLabuschagne:Survey07}. The reader can also find a nice exposition (especially, on the non-commutative Hilbert transform in the framework of $H^\infty(M,\tau)$) in Pisier and Xu's survey on non-commutative $L^p$-spaces \cite[\S8]{PisierXu:Handbook03}. 

More precisely, what we want to prove here is that for any non-zero singular $\varphi \in M^*$ in the sense of Takesaki \cite{Takesaki:TohokuMathJ58} one can find a ``peak" projection $p$ for $A$ in the sense of Hay \cite{Hay:IntegralEqOpTh07} such that $p$ dominates the (right) support projection of $\varphi$ but is smaller than the central support projection $z_s \in M^{\star\star}$ of the singular part $M^\star\ominus M_\star$. This is not exactly same as Amar and Lederer's result, but is enough for usual applications (even in classical theory for $H^\infty(\mathbb{D})$). Indeed, we will demonstrate it by proving that any non-commutative $H^\infty$-algebra $A = H^\infty(M,\tau)$ has the unique predual $M_\star/A_\perp$ with $A_\perp := \{\psi \in M_\star : \psi|_A = 0\}$. Proving it is our initial motivation; in fact, it can be regarded as an affirmative answer to the following natural (at least for us) question: Is the relative topology on $A$ induced from $\sigma(M,M_\star)$, which is most important, an intrinsic one of $A$ ? Also, our unique predual result may provide a new perspective in the direction of establishing the uniqueness of preduals by Grothendieck \cite{Grothendieck:CanadianJMath55} for $L^\infty$-spaces, by Dixmier \cite{Dixmier:BullFrance53} and Sakai \cite{Sakai:Pacific56} for von Neumann algebras or $W^*$-algebras, and then by Ando \cite{Ando:CommentMath78} and also a little bit later but independent work due to Wojtaszczyk \cite{Wojtaszczyk:Studia79} for $H^\infty(\mathbb{D})$. In particular, our result can be regarded as a simultaneous generalization of those classical results. Moreover, our result is an affirmative answer to a question posed by Godefroy stated in \cite{BlecherLabuschagne:Survey07}, and more importantly it covers any existing generalization like \cite{Chaumat:CRAcadParis79},\cite{Godefroy:TAMS84} of the above-mentioned work for $H^\infty(\mathbb{D})$ as a particular case. A natural ``Lebesgue decomposition" or ``normal/singular decomposition" for the dual of $H^\infty(M,\tau)$ is also given. The decomposition was first given by our ex-student Shintaro Sewatari in his master thesis \cite{Sewatari:MasterThesis} as a simple application of the non-commutative F.~and M.~Riesz theorem recently established by Blecher and Labuschagne \cite{BlecherLabuschagne:Studia07} so that the finite dimensionality assumption for the diagonal subalgebra $D$ was necessary there. Here it is established in full generality based on our Amar--Lederer type result instead of the non-commutative F.~and M.~Riesz theorem. After the completion of the presented work, the author found the paper \cite{Pfitzner:BLMS07} of H.~Pfitzner, where it is shown that any separable $L$-embedded Banach space $X$ becomes the unique predual of its dual $X^\star$. This means that establishing the Lebesgue decomposition is enough to show the uniqueness of predual for any non-commutative $H^\infty$-algebra $A = H^\infty(M,\tau)$ with  $M_\star$ separable. 

Our Amar--Lederer type result also enables us to remove the finite dimensionality assumption for the diagonal subalgebra $D$ from the results in \cite{BlecherLabuschagne:Studia07} numbered 3.5, 4.1, 4.2 and 4.3 there, including the non-commutative Gleason--Whitney theorem. Moreover, it gives a nice variant of Blecher and Labuschagne's non-commutative F.~and M.~Riesz theorem. Thus, it unexpectedly brings the current theory of non-commutative $H^p$-spaces due to Blecher and Labuschagne (see \cite{BlecherLabuschagne:Survey07}), which was already somewhat complete and satisfying, to an even more perfect and satisfactory form, though the presented work was initially aimed to prove the unique predual result for $H^\infty(M,\tau)$ as mentioned above. 
 
In closing, we should note that a bit different syntax has been (and will be) used for dual spaces. For a Banach space $X$ we denote by $X^\star$ and $X_\star$ its dual and predual instead of the usual $X^*$ and $X_*$, while $X^*$ stands for the set of adjoints of elements in $X$ when $X$ is a subset of a $C^*$-algebra.          

\medskip\noindent
{\it Acknowledgment.} We thank Professor Timur Oikhberg for kindly advising us to mention what the unique predual $M_\star/A_\perp$ possesses Pelczynski's property {\rm(V$^*$)} in Corollary \ref{C3.3} explicitly. We also thank the anonymous referee for his or her critical reading and a number of fruitful suggestions, which especially enable us to improve the presentation of the materials given in \S4.   

\section{Amar--Lederer Type Result for $H^\infty(M,\tau)$}

Let $A = H^\infty(M,\tau)$ be a non-commutative $H^\infty$-algebra with a finite von Neumann algebra $M$ and a faithful normal tracial state $\tau$ on $M$.

\begin{theorem} \label{T2.1} For any non-zero singular $\varphi \in M^\star$ there is a contraction $a \in A$ and a projection $p \in M^{\star\star}$ such that 
\begin{itemize}
\item[(2.1.1)] $a^n$ converges to $p$ in the $w^*$-topology $\sigma(M^{\star\star},M^\star)$ as $n\rightarrow\infty$; 
\item[(2.1.2)] $\langle|\varphi|,p\rangle = |\varphi|(1)$; 
\item[(2.1.3)] $\langle\psi,p\rangle = 0$ for all $\psi \in M_\star$ {\rm(}regarded as a subspace of $M^\star${\rm)}, or equivalently $a^n$ converges to $0$ in $\sigma(M,M_\star)$ as $n\rightarrow\infty$. This, in particular, shows that $p \le z_s$.  
\end{itemize}  
Here, $\langle\cdot,\cdot\rangle : M^\star \times M^{\star\star} \rightarrow \mathbf{C}$ is the dual pairing and $|\varphi|$ denotes the absolute value of $\varphi$ with the polar decomposition $\varphi = v\cdot|\varphi|$ due to Sakai \cite{Sakai:ProcJapanAcad58} and Tomita \cite{Tomita:MathJOkayama59}, when regarding $\varphi$ as an element in the predual of the enveloping von Neumann algebra $M^{\star\star}$ by $(M^{\star\star})_\star = M^\star$. 
\end{theorem}
\begin{proof}
Note that $|\varphi|$ is still singular. In fact,  $|\varphi| = v^*\cdot\varphi \in v^* z_s M^\star \subset z_s M^\star$ since $z_s$ is a central projection. Here $z_s$ stands for the central support projection of $M^\star\ominus M_\star$ as in \S1. The orthogonal families of non-zero projections in $\mathrm{Ker}|\varphi|$ clearly form an inductive set by inclusion, and then Zorn's lemma ensures the existence of a maximal family $\{q_k\}$, which is at most countable since $M$ is $\sigma$-finite. Let $q_0 := \sum_k q_k$ in $M$. If $q_0\neq1$, then Takesaki's criterion \cite{Takesaki:PJapanAcad59} shows the existence of a non-zero projection $r \in \mathrm{Ker}|\varphi|$ with $r \le 1-q_0$, a contradiction to the maximality. Thus, $q_0 = 1$. Moreover, if $\{q_k\}$ is a finite set, then $|\varphi|(1) = \sum_k |\varphi|(q_k) = 0$, a contradiction. Therefore, $\{q_k\}$ is a countably infinite family with $\sum_k q_k = 1$ in $M$. Letting $p_n := 1-\sum_{k\le n} q_k$ we have $p_n \rightarrow 0$ $\sigma$-weakly as $n\rightarrow\infty$ but $|\varphi|(p_n) = |\varphi|(1) $ for all $n$. Set $p_0 := \bigwedge_n p_n$ in $M^{\star\star}$. Then, $\langle|\varphi|,p_0\rangle = \lim_n \langle|\varphi|,p_n\rangle = \lim_n |\varphi|(p_n) = |\varphi|(1) \neq 0$, and in particular,  $p_0 \neq 0$. 

Choosing a subsequence if necessary, we may and do assume $\tau(p_n) \le n^{-6}$. Then we can define an element $g := \sum_{n=1}^\infty n p_n \in L^2(M,\tau)$, the non-commutative $L^2$-space associated with $(M,\tau)$, since $\sum_{n=1}^\infty \Vert n p_n \Vert_{2,\tau} \le \sum_{n=1}^\infty n^{-2} < +\infty$. By the non-commutative Riesz theorem \cite[Theorem 1]{Randrianantoanina:JAustral98} and \cite[Theorem 5.4]{MarsalliWest:JOT98} there is an element $\tilde{g}=\tilde{g}^* \in L^2(M,\tau)$ such that $f := g+\sqrt{-1}\tilde{g}$ falls in the closure $[A]_{2,\tau}$ of $A$ in $L^2(M,\tau)$ via the canonical embedding $M\hookrightarrow L^2(M,\tau)$. We can regard $g,\tilde{g},f \in L^2(M,\tau)$ as unbounded operators, affiliated with $M$, on the Hilbert space $\mathcal{H} := L^2(M,\tau)$ with a common core $\mathcal{D}$. Let $\xi \in \mathcal{D}$ be chosen arbitrary. Since $g \ge 0$ and $\tilde{g} = \tilde{g}^*$, one has $\Vert(1+f)\xi\Vert_{2,\tau}\Vert\xi\Vert_{2,\tau} \ge |((1+f)\xi|\xi)_\tau| = |(\xi|\xi)_\tau + (g\xi|\xi)_\tau + \sqrt{-1}(\tilde{g}\xi|\xi)_\tau| \ge \Vert\xi\Vert_{2,\tau}^2$ and similarly $\Vert(1+f)^*\xi\Vert_{2,\tau}\Vert\xi\Vert_{2,\tau} \ge \Vert\xi\Vert_{2,\tau}^2$, and hence $(1+f)^{-1} \in M$ exists and $\Vert(1+f)^{-1}\Vert_\infty \le 1$. Also, similarly one has $\Vert(1+f)\xi\Vert_{2,\tau}\Vert\xi\Vert_{2,\tau} \ge |((1+f)\xi|f\xi)_\tau| = |(\xi|f\xi)_\tau + (f\xi|f\xi)_\tau| = |(\xi|g\xi)_\tau - \sqrt{-1}(\xi|\tilde{g}\xi)_\tau + (f\xi|f\xi)_\tau| \ge \Vert f\xi\Vert_{2,\tau}^2$ so that $\Vert f\xi\Vert_{2,\tau} \le \Vert(1+f)\xi\Vert_{2,\tau}$ holds. Therefore, $\Vert f(1+f)^{-1}\zeta\Vert_{2,\tau} \le \Vert\zeta\Vert_{2,\tau}$ for all $\zeta \in \mathcal{H}$, and thus $b:=f(1+f)^{-1} \in M$ is a contraction.  

We will then prove that $b$ actually falls in $A$. First, recall the following standard but non-trivial fact: any bounded element in the closure $[A]_{p,\tau}$ of $A$ in $L^p(M,\tau)$, the non-commutative $L^p$-space, falls in $A$. In fact, let $x \in [A]_{p,\tau}$ be a bounded element, i.e., $x \in M$, and then there is a sequence $\{a_n\}$ in $A$ with $\Vert a_n - x\Vert_{p,\tau} \longrightarrow    
0$ as $n\rightarrow\infty$. For each $y \in A$ with $E(y)=0$ one has $\Vert a_n y - xy\Vert_{p,\tau}  \longrightarrow 0$ as $n\rightarrow\infty$ so that $\tau(xy) = \lim_n\tau(a_n y) = 0$ implying $x \in A$, where we use $A = \{ x \in M : \tau(xy) = 0\ \text{for all $y \in A$ with $E(y)=0$}\}$ due to Arveson \cite{Arveson:AJM67}. (It seems that this fact is used but not mentioned explicitly in the final step of the proof of \cite[Lemma 2]{Randrianantoanina:JAustral98} that we need here). Letting $g_N := \sum_{n=1}^N np_n \in M$ with its conjugate $\widetilde{g_N}$ we have $f_N = g_N + \sqrt{-1}\widetilde{g_N} \longrightarrow f$ in $L^2(M,\tau)$ as $N\rightarrow\infty$ thanks to the non-commutative Riesz theorem \cite[Theorem 1]{Randrianantoanina:JAustral98} and \cite[Theorem 5.4]{MarsalliWest:JOT98}. As before, for each $N$ one has $(1+f_N)^{-1} \in M$ and $\Vert(1+f_N)^{-1}\Vert_\infty\le1$, and moreover the discussion in \cite[Lemma 2]{Randrianantoanina:JAustral98} shows that $(1+f_N)^{-1}$ indeed falls in $A$. Since $(1+f)^{-1} \in M$ and $\Vert(1+f)^{-1}\Vert_\infty\le1$ as shown before, we have, for each $\xi \in M \subset L^2(M,\tau)$ (a right-bounded vector in $L^2(M,\tau)$), $\Vert((1+f_N)^{-1}-(1+f)^{-1})\xi\Vert_{2,\tau} = \Vert (1+f_N)^{-1} (f-f_N)(1+f)^{-1}\xi\Vert_{2,\tau} \le \Vert \xi \Vert_\infty \Vert f - f_N\Vert_{2,\tau} \longrightarrow 0$ as $N\rightarrow\infty$ so that $(1+f)^{-1} = \lim_N (1+f_N)^{-1} \in A$ in strong operator topology, implying $b = f(1+f)^{-1} \in M \cap [A]_{2,\tau} = A$ as claimed above.  

As before we have $\Vert(1+f)\xi\Vert_{2,\tau}\Vert\xi\Vert_{2,\tau} \ge |((1+f)\xi|\xi)_\tau| \ge (g\xi|\xi)_\tau \ge n(p_n\xi|\xi)_\tau = n\Vert p_n\xi\Vert_{2,\tau}^2$ for each $\xi \in \mathcal{D}$. Here the inequality $(g\eta|\eta)_\tau \ge n(p_n\eta|\eta)_\tau$ for $\eta$ in the domain of $g$ is used. (This can be easily checked when $\eta$ is in $M \subset L^2(M,\tau)$, and $M \subset L^2(M,\tau)$ is known to form a core of $g$ thanks to a classical result, see, e.g.~\cite[Theorem 9.8]{StratilaZsido:Book}). Thus, letting $\xi := (1+f)^{-1}\zeta$ for each $\zeta \in \mathcal{H}$ we get $\Vert p_n(1+f)^{-1}\zeta\Vert_{2,\tau}^2 \le n^{-1}\Vert\zeta\Vert_{2,\tau}\Vert(1+f)^{-1}\zeta\Vert_{2,\tau} \le n^{-1}\Vert\zeta\Vert_{2,\tau}^2$ so that $\Vert p_n - p_n b \Vert_\infty = \Vert p_n(1+f)^{-1}\Vert_\infty \le n^{-1/2}$. In the universal representation $M \curvearrowright \mathcal{H}_u$ we have $\Vert(p_0 - p_0 b)\zeta\Vert_{\mathcal{H}_u} \le \Vert p_0\zeta - p_n\zeta\Vert_{\mathcal{H}_u} + \Vert p_n - p_n b\Vert_\infty\Vert\zeta\Vert_{\mathcal{H}_u} + \Vert p_n(b\zeta) - p_0(b\zeta)\Vert_{\mathcal{H}_u} \le \Vert p_0\zeta - p_n\zeta\Vert_{\mathcal{H}_u} + n^{-1/2}\Vert\zeta\Vert_{\mathcal{H}_u} + \Vert p_n(b\zeta) - p_0(b\zeta)\Vert_{\mathcal{H}_u} \longrightarrow 0$ as $n \rightarrow \infty$ for each $\zeta \in \mathcal{H}_u$ since $p_0 = \bigwedge_n p_n$ in $M^{\star\star} = M''$ on $\mathcal{H}_u$. Since $b$ is a contraction, we get $p_0 = p_0 b = bp_0 = p_0 bp_0$. Then, by \cite[Lemma 3.7]{Hay:IntegralEqOpTh07} the new contraction $a := (1+b)/2$ satisfies that $a^n$ converges to a certain projection $p \in M^{\star\star}$ in $\sigma(M^{\star\star},M^\star)$ as $n\rightarrow\infty$, and $p_0 \le p$ so that $\langle|\varphi|,p\rangle = |\varphi|(1)$. If a vector $\xi \in \mathcal{H}$ satisfies $\Vert a\xi\Vert_{2,\tau} = \Vert\xi\Vert_{2,\tau}$, then $2\Vert\xi\Vert_{2,\tau} = \Vert \xi + b\xi\Vert_{2,\tau} \le \Vert\xi\Vert_{2,\tau}+\Vert b\xi\Vert_{2,\tau} \le 2\Vert\xi\Vert_{2,\tau}$, which implies $\Vert b\xi\Vert_{2,\tau} = \Vert\xi\Vert_{2,\tau}$ and $\Vert\xi + b\xi\Vert_{2,\tau} = \Vert\xi\Vert_{2,\tau}+\Vert b\xi\Vert_{2,\tau}$. Then, it is plain to see that these two norm conditions imply $b\xi=\xi$. However, $(1+f)^{-1}\xi = (1-b)\xi = 0$ so that $\xi = 0$. Therefore, there is no reducing subspace of $a$ in $\mathcal{H}$, on which $a$ acts as a unitary. Hence the so-called Foguel decomposition (\cite{Foguel:PacificJMath63}) shows that $a^n \longrightarrow 0$ $\sigma$-weakly as $n\rightarrow\infty$. In particular, $\langle \psi,p\rangle = \lim_n \langle\psi,a^n\rangle = \lim_n \psi(a^n) = 0$ for all $\psi \in M_\star$.    
\end{proof}

Choose $\varphi \in M^\star$, and decompose it into the normal and singular parts $\varphi = \varphi_n + \varphi_s$ with $\varphi_n := (1-z_s)\cdot\varphi \in M_\star$ and $\varphi_s := z_s\cdot\varphi \in M^\star\ominus M_\star$. Assume that $\varphi_s \neq 0$, and let $p \in M^{\star\star}$ be a projection for $\varphi_s$ as in Theorem \ref{T2.1}. By (2.1.2) and the polar decomposition $\varphi_s = v\cdot|\varphi_s|$ we have $|\langle \varphi_s,(1-p)x\rangle| =  |\langle v\cdot|\varphi_s|,(1-p)x\rangle| \le \langle|\varphi_s|,1-p\rangle^{1/2}\langle|\varphi_s|,v^* x^* x v\rangle^{1/2} = 0$ for every $x \in M^{\star\star}$ so that $\varphi_s\cdot(1-p) = 0$, i.e., $\varphi_s = \varphi_s\cdot p$. Moreover, by (2.1.3) a similar estimate shows $\varphi_n\cdot p = 0$. Hence, we get $\varphi_s = \varphi\cdot p$. Therefore we have the following corollary: 

\begin{corollary}\label{C2.2} If $\varphi \in M^\star$ has the non-zero singular part $\varphi_s \in M^\star\ominus M_\star$, then there is a contraction $a \in A$ and a projection $p \in M^{\star\star}$ such that $a^n \longrightarrow p$ in $\sigma(M^{\star\star},M^\star)$, $a^n \longrightarrow 0$ in $\sigma(M,M_\star)$  as $n\rightarrow\infty$ and $\varphi_s = \varphi\cdot p$.
\end{corollary}

We next examine the contraction $a$ and the projection $p$ in Theorem \ref{T2.1} and/or Corollary \ref{C2.2}. By \cite[Lemma 3.6]{Hay:IntegralEqOpTh07}, $a$ peaks at $p$ and moreover $(a^* a)^n \searrow p$ in $\sigma(M^{\star\star},M^\star)$ as $n\rightarrow\infty$ so that $p$ is a closed projection in the sense of Akemann \cite{Akemann:JFA69},\cite{Akemann:JFA70}. For any positive $\psi \in M^\star$ one has $\sum_{n=2}^N |\psi((a^* a)^n - (a^* a)^{n-1})| = -\sum_{n=2}^N \psi((a^* a)^n - (a^* a)^{n-1}) = \psi(a^* a) - \psi((a^* a)^N) \longrightarrow \langle \psi, a^* a - p\rangle$ as $N\rightarrow\infty$, from which one easily observes that the sequence $\{(a^* a)^n\}$ is weakly unconditionally convergent, see, e.g.~\cite[D\'{e}finition 1]{GodefroyTalagrand:CRParis81}. This fact is necessary in the course of proving that $M_\star/A_\perp$ is the unique predual of $A$.  
 
\section{First Applications: Predual of $H^\infty(M,\tau)$}

We first establish the following theorem: 

\begin{theorem}\label{T3.1} $M_\star/A_\perp$ is the unique predual of $A=H^\infty(M,\tau)$. 
\end{theorem}

A Banach space $E$ is said to have a unique predual when the following property holds: If the duals $F^\star$ and $G^\star$ of two other Banach spaces $F$ and $G$ are isometrically isomorphic to $E$, then $F = G$ must hold in the dual $E^{\star}$ via the canonical embeddings. Our discussion will be done in the line presented in \cite[IV]{Godefroy:TAMS84} so that what we will actually prove is that $M_\star/A_\perp$ has property (X) in the sense of Godefroy and Talagrand and the desired assertion immediately follows from their result, see \cite[D\'{e}finition 3, Th\'{e}or\`{e}me 5]{GodefroyTalagrand:CRParis81}. 

\begin{proof} Choose $\varphi \in A^\star$, and then one can extend it to $\tilde{\varphi} \in M^\star$ by the Hahn--Banach extension theorem. Decompose $\tilde{\varphi}$ into the normal and singular parts $\tilde{\varphi} = \tilde{\varphi}_n + \tilde{\varphi}_s$. It suffices to show the following: If $\lim_n \varphi(x_n) = 0$ for any weakly unconditionally convergent sequence $\{x_n\}$ in $A$ with $x_n \longrightarrow 0$ in $\sigma(A,M_\star/A_{\perp})$ or the relative topology from $\sigma(M,M_\star)$ as $n\rightarrow\infty$, then $\tilde{\varphi}_s|_A = 0$, that is, $\varphi = \tilde{\varphi}_n|_A$ must hold. We may assume $\tilde{\varphi}_s \neq 0$. By Corollary \ref{C2.2} together with the discussion just below it, we can find two sequences $\{a_n\}$ and $\{b_n\}$ and a projection $p \in M^{\star\star}$ such that (i) the $a_n$'s are in $A$; (ii) the $b_n$'s are in $M$ and $\{b_n\}$ is weakly (in $\sigma(M,M^\star)$) unconditionally convergent; (iii) both $a_n$ and $b_n$ converge to $p$ in $\sigma(M^{\star\star},M^\star)$ but to $0$ in $\sigma(M,M_\star)$; (iv) $\tilde{\varphi}_s = \tilde{\varphi}\cdot p$. Then, as same as in \cite[Th\'{e}or\`{e}me 33]{Godefroy:TAMS84} (by using a trick in \cite[the proof of Proposition 1.c.3 in p.32]{LindenstraussTzafriri:BookVolII}) we may and do assume that $\{a_n\}$ is also weakly unconditionally convergent. Let $x \in A$ be chosen arbitrary, and then $\{a_n x\}$ clearly becomes weakly unconditionally convergent. Moreover, it trivially holds that $a_n x \longrightarrow 0$ in $\sigma(M,M_\star)$ as $n\rightarrow\infty$. Therefore, we have $\tilde{\varphi}_s(x) = \langle \tilde{\varphi},px\rangle = \lim_n \langle \tilde{\varphi},a_n x\rangle = \lim_n \varphi(a_n x) = 0$ by the assumption here. 
\end{proof}  

As is well-known the predual $M_\star$ of a von Neumann algebra $M$ can be naturally embedded to the dual $M^\star$ as the range of an $L$-projection, see \cite{Takesaki:TohokuMathJ58}. Hence it is natural to ask whether the predual $M_\star/A_\perp$ of $A = H^\infty(M,\tau)$ can be also embedded to the dual $A^\star$ as the range of an $L$-projection. This is indeed true in general. Here we will explain it as an application of our Amar--Lederer type result. 

Denote by $A_n^\star$ the set of all $\varphi \in A^\star$ that can be extended to $\tilde{\varphi} \in M_\star$, and also by $A_s^\star$ the set of all $\psi \in A^\star$ that can be extended to $\tilde{\psi} \in M^\star\ominus M_\star$. This definition agrees with \cite[p.35]{Ando:CommentMath78}. For any $\varphi \in A^\star$, by the Hahn--Banach extension theorem one can extend it to $\tilde{\varphi} \in M^\star$. Then, decompose $\tilde{\varphi}$ into the normal and singular parts $\tilde{\varphi} = \tilde{\varphi}_n + \tilde{\varphi}_s$. We set $\varphi_n := \tilde{\varphi}_n|_A \in A_n^\star$ and $\varphi_s := \tilde{\varphi}_s|_A \in A_s^\star$. Then we call $\varphi = \varphi_n + \varphi_s$ an ``$(M\supset A)$-Lebesgue decomposition" of $\varphi$. On first glance, it is likely that this decomposition depends on the particular choice of the extension $\tilde{\varphi}$. However, we have:   

\begin{proposition}\label{P3.2} The following hold true{\rm:}
\begin{itemize}
\item[(3.4.1)] $A_n^\star\cap A_s^\star = \{0\}$. 
\item[(3.4.2)] The notion of $(M\supset A)$-Lebesgue decomposition $\varphi = \varphi_n+\varphi_s$ of $\varphi \in A^\star$ is well-defined, that is, $\varphi_n$ and $\varphi_s$ are uniquely determined by $\varphi$. Moreover, $\Vert\varphi\Vert = \Vert\varphi_n\Vert + \Vert\varphi_s\Vert$ holds. 
\end{itemize} 
\end{proposition}
\begin{proof} (3.4.1) On contrary, suppose that there is a non-zero $\varphi \in A_n^\star\cap A_s^\star$, and then one can choose $\tilde{\varphi}_n \in M_\star$ and $\tilde{\varphi}_s \in M^\star\ominus M_\star$ in such a way that $\varphi = \tilde{\varphi}_n|_A = \tilde{\varphi}_s|_A$. Since $\varphi \neq 0$ implies  $\tilde{\varphi}_s\neq0$, one can find, by Corollary \ref{C2.2}, a contraction $a \in A$ and a projection $p \in M^{\star\star}$ so that $a^n \longrightarrow p$ in $\sigma(M^{\star\star},M^\star)$, $a^n \longrightarrow 0$ in $\sigma(M,M_\star)$ as $n\rightarrow\infty$ and $\tilde{\varphi}_s = \tilde{\varphi}_s\cdot p$. Let $x \in A$ be arbitrary, and $a^n x \longrightarrow 0$ in $\sigma(M,M_\star)$ clearly holds. Then one has $\varphi(x) = \tilde{\varphi}_s(x) = \langle\tilde{\varphi}_s,px\rangle = \lim_n \langle\tilde{\varphi}_s,a^n x\rangle = \lim_n \varphi(a^n x) = \lim_n \tilde{\varphi}_n(a^n x) = 0$, a contradiction.   

(3.4.2) Assume that we have two $(M\supset A)$-Lebesgue decompositions $\varphi = \varphi_{n1} + \varphi_{s1} = \varphi_{n2} + \varphi_{s2}$. Then $\varphi_{n1}-\varphi_{n2} = \varphi_{s2}-\varphi_{s1} \in A_n^\star\cap A_s^\star = \{0\}$ by (3.4.1) so that  $\varphi_{n1}=\varphi_{n2}$ and $\varphi_{s1}=\varphi_{s2}$. Hence the $(M\supset A)$-Lebesgue decomposition is well-defined. Let $\tilde{\varphi} \in M^\star$ be the Hahn-Banach extension of $\varphi$, i.e., $\Vert\tilde{\varphi}\Vert = \Vert\varphi\Vert$. By definition we have $\varphi_n = \tilde{\varphi}_n|_A$ and $\varphi_s = \tilde{\varphi}_s|_A$. Then one has $\Vert\varphi\Vert = \Vert\tilde{\varphi}\Vert = \Vert\tilde{\varphi}_n\Vert + \Vert\tilde{\varphi}_s\Vert \ge \Vert\varphi_n\Vert + \Vert\varphi_s\Vert \ge \Vert \varphi_n + \varphi_s\Vert = \Vert\varphi\Vert$ so that the desired norm equation follows.      
\end{proof}
 
\begin{corollary}\label{C3.3} The predual $M_\star/A_\perp$ of $A = H^\infty(M,\tau)$ is the range of an $L$-projection from $A^\star$. Hence $M_\star/A_\perp$ has Pe\l czy\'nski's property {\rm(V$^*$)}, and, in particular, is sequentially weakly complete. 
\end{corollary} 
\begin{proof} The first part is immediate from the above proposition since $A_n^\star = M_\star/A_\perp$ trivially holds. The latter half is due to Pfitzner's theorem \cite{Pfitzner:Studia93} and an observation of Pe{\l}czy\'nski \cite[Proposition 6]{Pelczynski:BullAcadPolon62}. 
\end{proof}

It seems a natural question to find an ``intrinsic characterization" of singularity for elements in $A^\star$ like Takesaki's criterion \cite{Takesaki:PJapanAcad59}.  It seems that there is no such result even in the classical theory.  

\section{Second Applications: Noncommutative Function Algebra Theory} 

In this section we will explain how our Amar--Lederer type result nicely complements the non-commutative function algebra theory due to Blecher and Labuschagne \cite{BlecherLabuschagne:Studia07}. The key is the following variant of Blecher and Labuschagne's F.~and M.~Riesz theorem, which was given implicitly in the previous version of these notes. The current, explicit formulation was suggested by the referee.  

\begin{theorem}\label{T4.1} Any non-commutative $H^\infty$-algebra $A=H^\infty(M,\tau)$ satisfies the following property: Whenever $\varphi \in M^\star$ annihilates $A$, the normal and singular parts $\varphi_n$ and $\varphi_s$ annihilate $A$ separately. 
\end{theorem}  
\begin{proof} Although the proof is essentially same as that of Proposition \ref{P3.2}, we do give it for the sake of completeness. Let us choose $\varphi \in M^\star$ in such a way that $\varphi|_A = 0$, and decompose it into the normal and singular parts $\varphi = \varphi_n + \varphi_s$. On contrary, we assume that $\varphi_n|_A \neq0$ or $\varphi_s|_A\neq0$. If there existed $x \in A$ with $\varphi_n(x) \neq 0$, then it would follow that $\varphi_s(x) = -\varphi_n(x) \neq 0$. Thus we may assume that $\varphi_s|_A \neq 0$. Then, by Corollary \ref{C2.2} one can find a contraction $a \in A$ and a projection $p \in M^{\star\star}$ so that $a^n \longrightarrow p$ in $\sigma(M^{\star\star},M^\star)$, $a^n \longrightarrow 0$ in $\sigma(M,M_\star)$ as $n\rightarrow\infty$ and $\varphi_s = \varphi_s\cdot p$. For any $x \in A$, the $a^n x$'s still fall in $A$ but $\varphi(a^n x) \longrightarrow \langle\varphi,px\rangle = \varphi_s(x)$, and consequently $\varphi_s(x) = 0$, a contradiction.    
\end{proof}

Blecher and Labuschagne's F.~and M.~Riesz theorem \cite[\S3]{BlecherLabuschagne:Studia07}, which is originated in classical theory, asserts a quite similar property, that is, whenever $\varphi \in M^\star$ annihilates $A_0 := \{ a \in A : E(a) = 0\}$ the normal and singular parts $\varphi_n$ and $\varphi_s$ annihilate $A_0$ and $A$, respectively, and moreover its necessary and sufficient requirement is that $D$ is finite dimensional. Note that this F.~and M.~Riesz property is apparently stronger than the consequence of Theorem \ref{T4.1} here, and it should be remarked that the proofs of Corollary 3.5, Theorem 4.1, Theorem 4.2 and Corollary 4.3 in \cite{BlecherLabuschagne:Studia07} need only the consequence of the above theorem but do not use Blecher and Labuschagne's F.~and M.~Riesz theorem itself. Thus, they all hold true without any assumption. Consequently, we get the next theorem. 

\begin{theorem} Any non-commutative $H^\infty$-algebra $A=H^\infty(M,\tau)$ enjoys the following{\rm:} 
\begin{itemize}
\item[(4.2.1)] If $\varphi \in M^\star$ annihilates $A+A^*$, then the $\varphi$ must be singular. {\rm (}cf. \cite[Corollary 3.5]{BlecherLabuschagne:Studia07}.{\rm)}
\item[(4.2.2)] Every Hahn--Banach extension to $M$ of any normal {\rm (}i.e., continuous in the relative topology induced from $\sigma(M,M_\star)${\rm )} functional on $A$ must fall in $M_\star$. {\rm (}cf. the second part of \cite[Theorem 4.1]{BlecherLabuschagne:Studia07}. {\rm)}
\item[(4.2.3)] Any $\varphi \in M_\star$ is the unique Hahn--Banach extension of its restriction to $A+A^*$. In particular, $\Vert\varphi\Vert=\Vert\varphi|_{A+A^*}\Vert$ for any $\varphi \in M_\star$. {\rm (}cf. \cite[Theorem 4.2]{BlecherLabuschagne:Studia07}.{\rm)}
\item[(4.2.4)] Any element in $M$ can be $\sigma$-weakly approximated by a norm-bounded net consisting of elements in $A+A^*$. {\rm (}cf. \cite[Corollary 4.3]{BlecherLabuschagne:Studia07}.{\rm)} 
\end{itemize}
\end{theorem} 

The above (4.2.2), called the non-commutative Gleason--Whitney theorem, might sound a contradiction to what Pe\l czy\'nski pointed out in \cite[Proposition 6.3]{Pelczynski:CBMS76}, a comment to Amar and Lederer's result. However, this is not the case since the (4.2.2) mentions only Hahn--Banach extensions. 

\begin{remarks} Following the referee's suggestion let us call a subalgebra $A$ of a finite von Neumann algebra $M$ with a faithful normal tracial state an F.~and M.~Riesz algebra if it satisfies the consequence of Theorem \ref{T4.1} of these notes. (Clearly, $H^\infty(M,\tau)$ with $D$ finite dimensional has been an example of F.~and M.~Riesz algebra since the appearance of \cite{BlecherLabuschagne:Studia07}.) Then, any F.~and M.~Riesz algebra has (GW1) of \cite{BlecherLabuschagne:Studia07}, and furthermore does (GW) of \cite{BlecherLabuschagne:Studia07} if $A+A^*$ is $\sigma$-weak dense in $M$. Also, Corollary 3.5, Theorem 4.1, Theorem 4.2 and Corollary 4.3 in \cite{BlecherLabuschagne:Studia07} hold true if $A+A^*$ is $\sigma$-weak dense in $M$ too. The proofs in \cite{BlecherLabuschagne:Studia07} still work without any change. The referee communicated to us that he or she had noticed in 2007 these observations together with the fact that any F.~and M.~Riesz algebra with separable predual has unique predual. 
\end{remarks}

}

\end{document}